\input amstex
\documentstyle{amsppt}
\magnification=\magstep1

\NoBlackBoxes
\TagsAsMath

\pagewidth{6.5truein}
\pageheight{9.0truein}

\long\def\ignore#1\endignore{\par DIAGRAM\par}
\long\def\ignore#1\endignore{#1}

\ignore
\input xy \xyoption{matrix} \xyoption{arrow}
          \xyoption{curve}  \xyoption{frame}
\def\edge{\ar@{-}}
\def\dttdar{\ar@{.>}}

\def\dashedge{\ar@{--}}

\endignore

\def\la{{\Lambda}}
\def\lamod{\Lambda\text{-}\roman{mod}}
\def\Lamod{\Lambda\text{-}\roman{Mod}}
\def \len{\operatorname{length}}

\def\SS{{\Bbb S}}

\def\NN{{\Bbb N}}

\def\pdim{\operatorname{proj\,dim}}

\def\GL{\operatorname{GL}}

\def\Ext{\operatorname{Ext}}

\def\End{\operatorname{End}}

\def\findim{\operatorname{fin\, dim}}
\def\Findim{\operatorname{Fin\, dim}}
\def\pdim{\operatorname{p\, dim}}

\def\C{{\Cal C}}

\def\T{{\Cal T}}

\def\P{{\Cal P}}

\def\S{{\sigma}}

\def\T{{\Cal T}}

\def\V{{\frak V}}

\def\modlad{\operatorname{\bold{Mod}}_d(\Lambda)}
\def\modladelta{\operatorname{\bold{Mod}}_d(\Delta)}
\def\mod{\operatorname{\bold{Mod}}}

\def\pinflamod{\operatorname{\P^{< \infty}(\lamod)}}

\def\laySS{\operatorname{\bold{Mod}}(\SS)}
\def\lay#1{\operatorname{\bold{Mod}}(#1)}

\def\BHT{{\bf 1}}
\def\But{{\bf 2}}
\def\CB{{\bf 3}}
\def\CBS{{\bf 4}}
\def\Del{{\bf 5}}
\def\GKK{{\bf 6}}
\def\pre{{\bf 7}}
\def\dom{{\bf 8}}
\def\Sma{{\bf 9}}

\topmatter

\title Truncated path algebras are homologically transparent 
\endtitle

\rightheadtext{Generic representation theory}

\author A. Dugas, B. Huisgen-Zimmermann, and J. Learned
\endauthor

\address Department of Mathematics, University of California, Santa Barbara,
CA 93106-3080 \endaddress

\thanks The  second author was partly supported by a grant from the National
Science Foundation. 
\endthanks

\dedicatory  To the memory of Tony Corner
\enddedicatory

\abstract It is shown that path algebras modulo relations of the form $\Lambda = KQ/I$, where $Q$ is a quiver, $K$ a coefficient field, and $I \subseteq KQ$ the ideal generated by all paths of a given length, can be readily analyzed homologically, while displaying a wealth of phenomena. In particular, the syzygies of their modules, and hence their finitistic dimensions, allow for smooth descriptions in terms of $Q$ and the Loewy length of $\Lambda$. The same is true for the distributions of projective dimensions attained on the irreducible components of the standard parametrizing varieties for the modules of fixed $K$- dimension. 
\endabstract

\endtopmatter

\document

\head 1. Introduction and notation \endhead

The problem of opening up general access roads to the finitistic dimensions
of a finite dimensional algebra $\la$, given through quiver and relations, is
quite challenging.  This is witnessed, for instance, by the fact that the
longstanding question ``Is  the (left) little finitistic dimension of $\la$,
$$\findim \la = \sup\{\pdim M \mid M \in \pinflamod\},$$ always finite?"
(Bass 1960) has still not been settled.  Here
$\pdim M$ is the projective dimension of a module $M$, and $\pinflamod$
denotes the category of finitely generated (left) $\la$-modules of finite
projective dimension.

In \cite{\BHT}, Babson, the second author, and Thomas showed that truncated
path algebras of quivers are particularly amenable to geometric exploration,
while nonetheless displaying a wide range of interesting phenomena.  This
led the authors of the present paper to the serendipitous discovery that the
same is true for the homology of such algebras.  By a truncated path algebra
we mean an algebra of the form $KQ/I$, where $KQ$ is the path algebra of a
quiver $Q$ with coefficients in a field $K$ and $I
\subseteq KQ$ the ideal generated by all paths of a fixed length $L+1$.  In
particular, truncated path algebras are monomial algebras. In this case, the
finitistic dimensions are known to be finite (see \cite{\GKK}).  Our goal
here is to show how much more is true in the truncated scenario. 

Roughly, our three main results (Theorems 5, 11, and 15) show the following
for a truncated path algebra $\la$:
\smallskip

\noindent $\bullet$  The little and big finitistic dimensions of $\la$
coincide and can be determined through a straightforward computation from
$Q$ and $L$.  Moreover, from a minimal amount of structural data for a
$\la$-module $M$, namely the radical layering $\SS(M) = \bigl(J^l M / J^{l+1}
M \bigr)_{0 \le l \le L}$  (or, alternatively, any ``skeleton" of $M$), one
can determine the syzygies and projective dimension of $M$ in a purely
combinatorial way.   (See Theorems 2, 5, and the first part of Theorem 11
for finer information.)
\smallskip

\noindent $\bullet$  The ``generic projective dimension" of any irreducible
component $\C$ of one of the classical module varieties (see beginning of
Section 3) is readily obtainable from graph-theoretic data as well.  So is
the full spectrum of values of the function $\pdim$ attained on the class of
modules parametrized by $\C$.  In particular, it turns out that the supremum
of the finite values among the generic finitistic dimensions of the various
irreducible components equals $\findim \la$.  (See Theorems 11 and 15 for
detail.)  
\smallskip

The picture emerging from the main theorems will be supplemented in a
sequel, where it will be shown that the category $\pinflamod$ is
contravariantly finite in the full category of finitely generated
$\la$-modules, whenever $\la$ is a  truncated path algebra. 
\bigskip

We fix a positive integer $L$.  Throughout, $\la$ denotes a truncated path
algebra of Loewy length $L+1$, that is, $\la = KQ/ I$, where $K$ is a field,
$Q$ a quiver, and $I$ the ideal generated by all paths of length
$L+1$.  The Jacobson radical $J$ of $\la$ satisfies $J^{L+1} = 0$ by
construction.  A ({\it nonzero\/}) {\it path in
$\la$\/} is the $I$-residue of a path in $KQ\setminus I$, that is, the
$I$-residue of a path $p$ in $KQ$ of length at most $L$; so, in particular,
any path in $\la$ is a {\it nonzero\/} element of $\la$ under this
convention.  Clearly, the paths in $\la$ form a $K$-basis for $\la$.  Due to
the fact that $I$ is homogeneous with respect to the path-length grading of
$KQ$, defining the {\it length\/} of such a path
$p + I$ to be that of $p$, yields an unambiguous concept of length for the
elements of this basis.  A distinguished role is played by the paths $e_1,
\dots, e_n$ of length zero in $\la$:  They constitute a full set of
orthogonal primitive idempotents, which is in obvious one-to-one
correspondence with the vertices of $Q$.  We will identify each $e_i$ with
the corresponding vertex, and whenever we refer to a primitive idempotent in
$\la$, we will mean one of the $e_i$.  Then the left ideals $\la e_i$ and
their radical factors $S_i = \la e_i/Je_i$, for $1 \le i \le n$, constitute
full sets of isomorphism representatives for the indecomposable projective
and simple left $\la$-modules, respectively. 

Finally, we say that a path $p$ in $\la$ or in $KQ$ is an {\it initial
subpath\/} of a path $q$ if there is a path $p'$ with $q = p' p$; here the
product $p' p$ stands for ``$p'$ after $p$."

\head 2.  Accessing the standard homological dimensions of $\Lamod$ \endhead
The (left) {\it  big finitistic dimension\/} of $\la$ is the supremum,
$\Findim \la$, of the projective dimensions of all left $\la$-modules of
finite projective dimension; for the {\it little finitistic dimension\/}
consult the introduction.    We start by recording some prerequisites
established in
\cite{\BHT}.  As was shown in \cite{\BHT}, the well-known fact that all {\it
second\/} syzygies of modules over a monomial algebra are direct sums of
cyclic modules generated by paths of positive length (see \cite{\pre} and
\cite{\But}), can be improved for truncated path algebras so as to cover
{\it first\/} syzygies as well.  In particular, this makes the big and
little finitistic dimensions of $\la$ computable from a finite set of cyclic
test modules.   

More sharply: Given any left $\la$-module $M$, we can explicitly pin down a
decomposition of the syzygy $\Omega^1(M)$ into cyclics.  This description of
$\Omega^1(M)$ relies on a {\it skeleton of
$M$\/}.  Roughly speaking, this is a path basis for $M$ with the property
that the path lengths respect the radical layering, $(J^lM / J^{l+1} M)_{0
\le l \le L}$.  The concept of a skeleton, defined in
\cite{\BHT} in full generality, can be significantly simplified for a
truncated path algebra $\la$.

\definition{Definition 1: Skeleton of a $\la$-module $M$}  Fix  a projective
cover $P$ of $M$, say $P = \bigoplus_{r \in R} \la z_r$, where each
$z_r$ is one of the primitive idempotents in $\{e_1, \dots,
e_n\}$, tagged with a place number $r$ (the index set $R$ may be infinite). 
A {\it path of length\/} $l$ in $P$ is any element $p z_r
\in P$, where $p$ is a path of length $l$ in $\la$ which starts in $z_r$ (in
particular, the paths in $P$ are again nonzero).  Identify
$M$ with an isomorphic factor module of $P$, say $M = P/C$.

(a)  A {\it skeleton of\/} $M = P/C$ is a set $\sigma$ of paths in $P$ such
that for each $l \le L$, the residue classes $q  + J^l M$ of the paths $q$
of length $l$ in $\sigma$ form a $K$-basis for $J^l M/ J^{l+1} M$. 
Moreover, we require, that
$\sigma$ be closed under initial subpaths, that is, if $q = p'pz_r \in
\sigma$, then $pz_r$ in $\sigma$.

(b) A path $q$ in $P \setminus \sigma$ is called
$\sigma$-{\it critical\/} if it is of the form
$q = \alpha p z_r$, where $\alpha$ is an arrow and $pz_r$ a path in
$\sigma$.   
\enddefinition

In particular, the definition entails that, for any skeleton $\S$ of $M =
P/C$, the full set of residue classes $\{q + C \mid q \in \sigma\}$ forms a
basis for
$M$.  Furthermore, it is easily checked that every $\la$-module $M$ has at
least one skeleton, and only finitely many when $M$ is finitely generated 
(as long as we keep the projective cover $P$ fixed). 

\proclaim{Theorem 2. Known Facts} \cite{\BHT, Lemma 5.10}  If $M$ is any
nonzero left
$\la$-module with skeleton $\S$, then 
$$\Omega^1(M) \cong \bigoplus_{q\ \S\text{-}\text{critical}} \la q.$$ In
particular, $\Omega^1(M)$ is isomorphic to a direct sum of cyclic left
ideals generated by nonzero paths of positive length in $\la$.  
\smallskip

\noindent Consequently, $\Findim \la = \findim \la = s + 1$, where 
$$s = \max\{ \pdim \la q \mid q \text{\ a path of positive length in
\ } \la\ \text{with}\ \pdim \la q < \infty \},$$  provided that the displayed
set is nonempty, and $s = -1$ otherwise. \qed
\endproclaim

We illustrate this result with an example which will accompany us
throughout. 

\definition{Example 3} Let $\la = KQ/I$, be the truncated path algebra of
Loewy length $L+1 = 4$ based on the following quiver $Q$

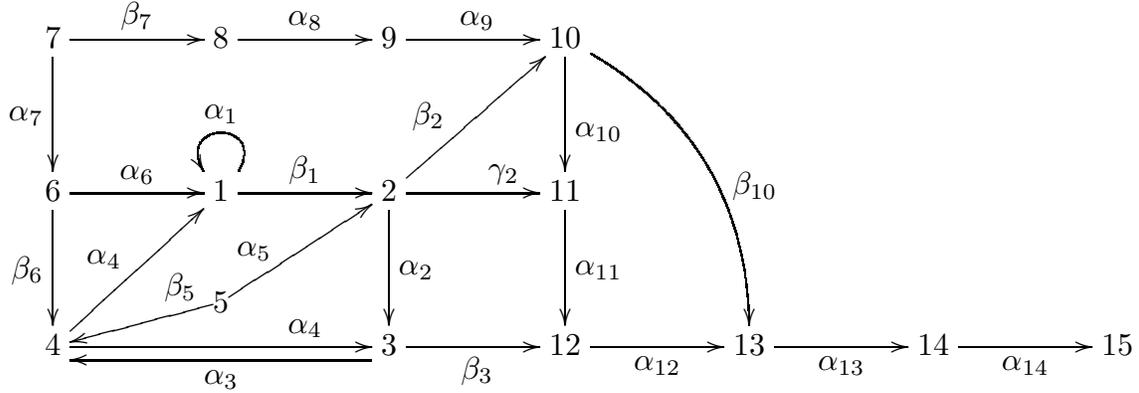
$$\xymatrixcolsep{3.5pc}\xymatrixrowsep{3pc}
\xymatrix{
7 \ar[d]_{\alpha_7} \ar[r]^{\beta_7} & 8 \ar[r]^{\alpha_8} & 9
\ar[r]^{\alpha_9} & 10 \ar[d]^(0.6){\alpha_{10}}
\ar@/^1.5pc/[ddr]^(0.6){\beta_{10}} \\ 
6 \ar[r]^{\alpha_6} \ar[d]_{\beta_6} & 1 \ar[r]^{\beta_1}
\ar@(ur,ul)[]_{\alpha_1} \save+<0ex,-8ex> \drop{5} \ar[r]^(0.3){\alpha_5}
\ar[dl]_(0.2){\beta_5} \restore & 2 \ar[r]^(0.65){\gamma_2}
\ar[ur]^(0.35){\beta_2} \ar[d]^{\alpha_2} & 11 \ar[d]^{\alpha_{11}} \\
4 \ar[rr]^(0.75){\alpha_4} \ar[ur]^(0.45){\alpha_4} && 3
\ar[ll]<1.0ex>^{\alpha_3} \ar[r]_{\beta_3} & 12 \ar[r]_{\alpha_{12}} & 13
\ar[r]_{\alpha_{13}} & 14 \ar[r]_{\alpha_{14}} & 15
}$$

\noindent Then the indecomposable projective left $\la$-modules $\la e_1$ and
$\la e_3$ have the following layered and labeled graphs (in the sense of
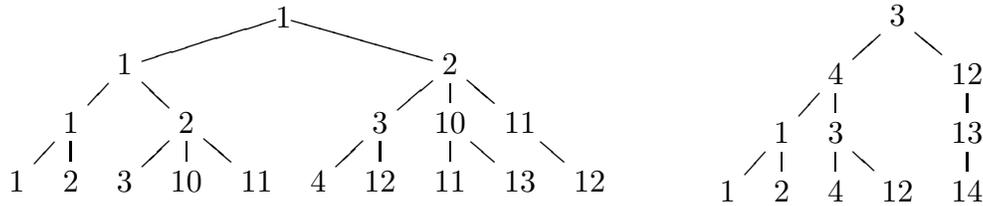
\cite{\pre} and \cite{\dom}):

$$\xymatrixrowsep{0.5pc} \xymatrixcolsep{0.5pc} \xymatrix{ & & & &
\save+<2ex,0ex> \drop{1} \ar@{-}[dll] \ar@{-}[drrr] \restore \\ & & 1
\ar@{-}[dl] \ar@{-}[dr] & & & & & 2 \ar@{-}[dl] \ar@{-}[d] \ar@{-}[dr] \\ &
1 \ar@{-}[dl] \ar@{-}[d] & & 2 \ar@{-}[dl] \ar@{-}[d] \ar@{-}[dr] &  & & 3
\ar@{-}[dl] \ar@{-}[d] & 10 \ar@{-}[d] \ar@{-}[dr] & 11 \ar@{-}[dr]
\\ 1 & 2 & 3 & 10 & 11 & 4 & 12 & 11 & 13 & 12} \ \ \  \ \ \  \ \ \ 
\xymatrixrowsep{0.5pc} \xymatrixcolsep{0.5pc} \xymatrix{ & & & 3 \ar@{-}[dl]
\ar@{-}[dr] \\ & & 4 \ar@{-}[dl] \ar@{-}[d] & & 12 \ar@{-}[d] \\ & 1
\ar@{-}[dl] \ar@{-}[d] & 3 \ar@{-}[d] \ar@{-}[dr] & & 13 \ar@{-}[d] \\ 1 & 2
& 4 & 12 & 14}$$

\noindent If $P = \la z_1$ with $z_1 = e_i$ in the notation of Definition 1,
each of the modules $\la e_i$ has a unique skeleton, which can be read off
the graph: It is the set of all initial subpaths of the edge paths in the
graph, read from top to bottom.  The skeleton of $\la e_1$, for instance,
consists of the paths $z_1 = e_1$ of length zero in
$P$, the paths 
$\alpha_1 z_1, \beta_1 z_1$ of length $1$, the paths  $\alpha_1^2 z_1,
\beta_1 \alpha_1 z_1, \alpha_2 \beta_1 z_1, \gamma_2 \beta_1 z_1, \beta_2
\beta_1 z_1$ of length $2$, together with all edge paths of length $3$.

For a sample application of Theorem 2, we consider the module $M$ determined
by the following graph:


$$\xymatrixrowsep{0.5pc} \xymatrixcolsep{0.5pc} 
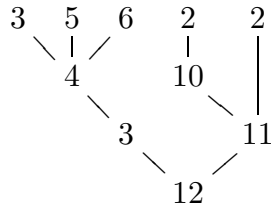
\xymatrix{3 \ar@{-}[dr] & 5
\ar@{-}[d] & 6 \ar@{-}[dl] &  2 \ar@{-}[d] & 2 \ar@{-}[dd]  \\ & 4
\ar@{-}[dr] & &  10 \ar@{-}[dr] \\ & & 3 \ar@{-}[dr] & & 11 \ar@{-}[dl] \\ &
& & 12}$$

\noindent  A projective cover of $M$ is $P = \la e_3 \oplus \la e_5 \oplus
\la e_6 \oplus (\la e_2)^2$, where $z_1 = e_3, z_2 = e_5$ and so on.  A
skeleton $\sigma$ of
$M$ (in this case there are several), together with the
$\sigma$-critical paths is communicated by the following graph, in which the
solid and dashed edges play different roles, as explained below: 


$$\xymatrixrowsep{1.5pc} \xymatrixcolsep{0.5pc} 
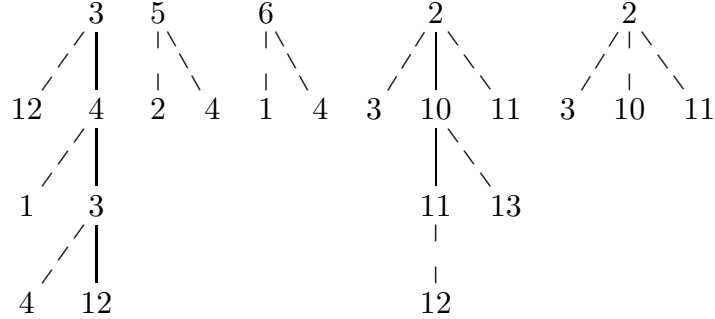
\xymatrix{ & 3 \ar@{-}[d]
\ar@{--}[dl] & 5 \ar@{--}[d] \ar@{--}[dr] & &  6 \ar@{--}[d] \ar@{--}[dr] &
& & 2 \ar@{--}[dl] \ar@{--}[dr] \ar@{-}[d] & & & 2 \ar@{--}[dl] \ar@{--}[dr]
\ar@{--}[d] \\ 12 & 4  \ar@{-}[d] \ar@{--}[dl] & 2 & 4 & 1 & 4 & 3 & 10
\ar@{-}[d] \ar@{--}[dr] & 11 & 3 & 10 & 11 \\ 1 & 3 \ar@{-}[d] \ar@{--}[dl]
& & & & & & 11  \ar@{--}[d] & 13 \\  4 & 12 & & & & & & 12}$$

\noindent  As above, the paths in $\sigma$ correspond to the intial subpaths
of the solidly drawn edge paths, including all paths of length zero -- e.g.,
$\beta_4 \alpha_3 z_1$, $z_3$ and $\alpha_{10} \beta_2 z_4$.  The
$\sigma$-critical paths are all the paths in the graph (again read from top
to bottom) which terminate in a dashed edge; for instance, $\alpha_3 \beta_4
\alpha_3 z_1$ and $\alpha_5 z_2$ are $\S$-critical.  Since
$\Omega^1(M) \cong \bigoplus_{q\ \S\text{-}\text{critical}} \la q$, we find
this syzygy to be the direct sum $\la \beta_3 \oplus \la \alpha_4 \alpha_3
\oplus \la \alpha_3 \beta_4 \alpha_3 \oplus \la \alpha_5 \oplus \la \beta_5
\oplus \la \alpha_6 \oplus \la \beta_6 \oplus \la \alpha _2 \oplus \cdots$. 
The graphs of $\la \alpha_4 \alpha_3$, $\la \beta_3$, and $\la \alpha_2$ are
respectively

$$\xymatrixrowsep{0.5pc} \xymatrixcolsep{0.5pc} 
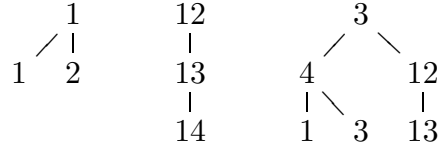
\xymatrix{& 1 \ar@{-}[dl]
\ar@{-}[d] \\ 1 & 2}  \ \ \ \ \ \ \  \xymatrixrowsep{0.5pc}
\xymatrixcolsep{0.5pc} \xymatrix{ 12 \ar@{-}[d] \\ 13 \ar@{-}[d] \\ 14} \ \
\ \ \ \ \ \xymatrixrowsep{0.5pc} \xymatrixcolsep{0.5pc} \xymatrix{ & 3
\ar@{-}[dl] \ar@{-}[dr] \\  4 \ar@{-}[dr] \ar@{-}[d] &  & 12 \ar@{-}[d] \\ 1
& 3  & 13}$$

\enddefinition

The main result of this section provides the projective dimensions of the
building blocks for the syzygies of arbitrary $\la$-modules; compare with
Theorem 2. 

\definition{Definition 4}  Let $l$ be a nonnegative integer $\le L$, and $c$
any nonnegative integer.  We define
$$l\text{-}\deg  (c) =  \left[\frac{c}{L+1}\right] +
\left[\frac{c+l}{L+1}\right].$$ Here $[x]$ stands for the largest integer
smaller than or equal to $x$.  Moreover, we set $l$-$\deg(\infty) = \infty$.
\enddefinition

The $l$-degree defines a nondecreasing function $\NN
\cup\{0, \infty\} \rightarrow \NN\cup\{0, \infty \}$ for any $l \le L$. 
Moreover, for $0 \le l \le l' \le L$ and arbitrary $c \in \NN \cup \{0\}$,
the difference $l'$-$\deg(c) - l$-$\deg(c) $ belongs to the set $\{0,1\}$. 
This observation will entail the final claim of the upcoming theorem, once
the first  --  displayed  --  equality is established. 

\proclaim{Theorem 5}  Suppose $q \in \la$ is a path of length $l > 0$ in
$\la$ {\rm(}i.e., the $I$-residue of a path of length at most $L$ in
$KQ${\rm )} with terminal vertex $e$.  Let $c = c(e)$ be the supremum of the
lengths of the paths in $KQ$ starting in $e$.  Then
$$\pdim \la q\   =   \ l\text{-}\deg(c).$$ In particular, $\pdim \la q <
\infty$ if and only if $c(e) < \infty$ {\rm (}meaning that there is no path
starting in $e$ and terminating on an oriented cycle{\rm )}.

Moreover, if $q'$ is another path in $\la$ that ends in $e$ such that $L \ge
\len(q') \ge \len(q) \ge 1$, then
$$\pdim \la q \le \pdim \la q' \le  1 +  \pdim \la q.$$
\endproclaim

In the Example, $c(e_7)$ is infinite, for instance, while $c(e_{10}) = 5$;
the latter shows that $\pdim (\la \alpha_9 \alpha_8 \beta_7) = 3$-$\deg(5) =
3$.  The argument backing Theorem 5 is purely combinatorial, the intuitive
underpinnings being of a graphical nature.  We start with two definitions
setting the stage. The first is clearly motivated by the statement of Theorem
5.   

\definition{Definition 6}  We call a vertex $e$ of the quiver $Q$ (alias a
primitive idempotent of $\la$) {\it cyclebound\/} in case there is a path
from $e$ to a vertex lying on an oriented cycle.  In case $e$ is cyclebound,
we also call the  simple module $\la e/ Je$ cyclebound.
\enddefinition  

Next, we consider the following partial order on the set of paths in
$KQ$.  Namely, given paths $p$ and $p'$ in $KQ$, we define
$$p' \le p\ \ \  \iff \ \ \ p' \text{\ is an initial subpath of\ } p;$$
recall that the latter amounts to the existence of a path $p''$ with the
property that $p = p'' p'$.  Hence, any two paths which are comparable have
the same starting point, and $e \le p$ for any path
$p$ starting in the vertex $e$.  Clearly, this partial order induces a
partial order on the set of paths in $\la$. 

Finally, we introduce a class of modules, which will turn out to tell the
full homological story of
$\la$.  The left ideals of the form $\la q$  --  the basic building blocks
of all syzygies of $\la$-modules  --  are among them.    

\definition{Definition 7 and comments.  Tree modules and branches} Any module
$\T$ of the form $\T \cong \la e/ V$, where $e$ is a vertex of $Q$ and $V =
(\sum_{v \in \V} \la v)$ is generated by some set $\V$ of paths of positive
length in
$\la e$ (possibly empty), will be called a {\it tree module with root\/}
$e$.  In particular, $\la e$ is a tree module with root $e$, the unique
candidate of maximal dimension among the tree modules with root $e$, in
fact; the simple module $\la e/Je$ is the tree module with root $e$ that has
minimal dimension.  

The terminology is motivated by the fact that the {\it graphs\/} of tree
modules are trees ``growing downwards" from their roots.    Note that tree
modules are determined up to isomorphism by their graphs. 

Given a tree module $\T$ as above, let $b_1,
\dots, b_r \in \la$ be the maximal paths in $\la e$  --  in the above partial
order  --  which are not contained in $V$.    The
$b_i$ are uniquely determined by the isomorphism class of
$M$ and are called the {\it branches\/} of $\T$.   Conversely, if we know $M$
to be a tree module, then the branches of $\T$ pin $\T$ down up to
isomorphism.  

If $\T \cong
\la e/Je$ is the simple tree module with root $e$, then $e$ is the only
branch of $\T$.  By contrast, if $\T = \la e/V$ is a nonsimple tree module,
then all branches of $\T$ have positive length.  Moreover, it is
straightforward to see that  $\T$ has a basis of the following form: 
$$\{e + V\} \cup \{q + V \mid q \text{\ is an initial subpath of positive
length of one of\ } b_1, \dots b_r\},$$  where $b_1, \dots, b_r$ are the
branches of $\T$.  If we pull back this basis to a set of paths in the
projective cover
$\la e$ of $\T$, then $\sigma$ is a skeleton of
$\T$ in the sense of Definition 1 (the only one).
\enddefinition       

Apart from $M$, all the modules displayed in Example 3 are tree modules. 
Their branches are precisely the maximal edge paths in their graphs, read
from top to bottom.  The proof of the next lemma is straightforward and we
leave it to the reader.  

\proclaim{Lemma 8}  Whenever $q$ is a path in $\la$ ending in $e$, not
necessarily of positive length, the cyclic left ideal
$\la q$ is a tree module with root $e$.  More precisely: If $l = \len(q)$,
let $b_1,
\dots, b_r$ be the maximal candidates among the paths of length $\le L -l$
starting in
$e$.  Then $\la q = \la e/ V$, where 
$$V = \Omega^1(\la q) = \bigoplus_{\beta \text{\ an arrow},\, i \le r} \la
\beta b_i\, ,$$ and the $b_i$ are the branches of $\la q$.

In particular, if $l > 0$, then $\pdim \la q < \infty$ if and only if $e$ is
non-cyclebound. 
\qed
\endproclaim

Combined with Theorem 2, Lemma 8 shows that all syzygies of
$\la$-modules are direct sums of tree modules.  Contrasting the final
statement for $l >0$, we see that, for the path
$q = e$ of length zero, $\la q = \la e$ is projective, irrespective of the
positioning of $e$ in $Q$.  As for the other extreme: By Lemma 8, the simple
module 
$S = \la e / Je$ has infinite projective dimension precisely when it is
cyclebound.  In Example 3, the vertices $e_1,
\dots, e_7$ are cyclebound, while $e_8, \dots, e_{15}$ are not.  Hence
$S_1, \dots, S_7$ are precisely the simple modules of infinite projective
dimension.

Note that the only potential branches $b_i$ of length $< L - l$ of a tree
module $\la q$ as in Lemma 8 end in a sink of the quiver $Q$.

\demo{Proof of Theorem 5} As in the statement of the theorem, let
$q$ be a path of positive length $l \le L$ in $\la$, which ends in the
vertex $e$.  In light of the remark preceding Theorem 5, we only need to
show the equality
$\pdim \la q = l$-$\deg(c)$, where $c = c(e)$ is the supremum of the lengths
of the paths in $KQ$ starting in $e$.  If $e$ is cyclebound, this equality
follows from Lemma 8.  So let us assume that
$e$ is non-cyclebound  --  meaning $c < \infty$  --  and induct on $c$.  If
$c \le  L - l$, all of the branches of the tree module $\la q$ end in sinks
of the quiver $Q$.  We infer that $\la q \cong \la e$ in that case, whence
$\pdim \la q = 0 = l$-$\deg(c)$.

Now suppose $c > L-l$, and assume that $\pdim \la p' = l'$-$\deg(c(e'))$ for
all paths $p'$ of length $l' \le L$ in $\la$ that end in a non-cyclebound
vertex $e'$ of $Q$ with $c(e') < c$.  Using the notation of Lemma 8, we
obtain
$\Omega^1(\la q) = \bigoplus_{\beta, \,i \le r} \la \beta b_i$, where the
$b_i$ are the branches of the tree module $\la q$ and the
$\beta$ are arrows.  Since the lengths of the $b_i$ are bounded from above
by $L - l \le L - 1$, the paths in $KQ$ of the form $\beta b_i$ where
$\beta$ is an arrow, have length at most $L$; therefore each of them gives
rise to a path in $\la$.  By the definition of $c$, there exists a path
$u$ of length $c$ in $KQ$ which starts in the vertex $e$, and by the
definition of the branches of $\la q$, there exists an index $j$ such that
$b_j$ is an initial subpath of $u$.   Necessarily,
$\len(b_j) = L-l$, because $\len(u) > L-l$.  In fact, $c > L-l$  guarantees
that  $u = u' \beta_j b_j$ in
$KQ$ for some arrow $\beta_j$ and a suitable path $u'$ of length $c' =
\len(u) - (L-l)  - 1 = c - (L-l)  - 1 \le c - 1$.   Since $u$ starts in the
non-cyclebound vertex $e$, the terminal vertex of $\beta_j b_j$  --  call it
$e'$  --  is again non-cyclebound.  Moreover, the maximality property of $u$
entails that $c' = c(e')$ is the maximal length of a path in $KQ$ starting
in $e'$.  Therefore, our induction hypothesis guarantees that $\pdim \la
\beta_j b_j = (L-l+1)\text{-}\deg  (c')$.  This degree in turn equals
$$\left[\frac{c'}{L+1}\right] + \left[\frac{c' + L - l + 1}{L+1}\right] =
\left[\frac{c + l - (L+1)}{L+1}\right] + \left[\frac{c}{L+1}\right] =
l\text{-}\deg(c) - 1;$$  the final equality follows from $\frac{c + l -
(L+1)}{L+1} = \frac{c+l}{L+1} - 1$.  Analogous applications of the induction
hypothesis, combined with the basic properties of the degree function, yield
$\pdim \la \beta b_i \le \pdim \la \beta_j b_j$ for any path $\beta b_i$
appearing in the decomposition of $\Omega^1(\la q)$.  We conclude that
$\pdim \la q = 1 + \pdim \la \beta_j b_j = l\text{-}\deg(c)$ as required. 
\qed
\enddemo

The following dichotomy for the finitistic dimension of $\la$ results from a
combination of Theorems 2 and 5 with Lemma 8.  
 
\proclaim{Corollary 9}  Suppose that $S_1, \dots, S_t$ are precisely the
non-cyclebound simple left $\la$-modules.  Then either 
$$\findim \la = \max_{1 \le i \le t} \pdim S_i \ \ \text{or}\ \  \findim
\la = 1 \, +\,\max_{1 \le i \le t} \pdim S_i,$$ and  
$$\max_{1 \le i \le m} \pdim S_i = 1\, +\, 1\text{-}\deg(m - 1),$$ where $m$
is the maximum of the lengths of the paths in $KQ$ which are not contingent
to any cycle.
\qed
\endproclaim

\noindent Both options for $\findim \la$ occur in concrete instances (see
below); of course, the smaller value equals the global dimension whenever
the quiver $Q$ is acyclic.  For the decision process in specific instances,
combine Theorems 2 and 5.  To contrast Corollary 9 with the homology of more
general algebras:  Recall that arbitrary natural numbers occur as finitistic
dimensions of monomial algebras all of whose simple modules have infinite
projective dimension.   So the corollary again attests to the degree of
simplification that occurs when the paths factored out of
$KQ$ have uniform length.   

\definition{Example 3 revisited} With the aid of Corollary 9, the finitistic
dimension of $\la$ can, in a first step, be computed up to an error of $1$,
through a simple count.  Here $m = 7$, and $L = 3$, whence the maximum of
the projective dimensions of the non-cyclebound simple modules (here $S_8,
\dots, S_{15}$) is $1 + 1$-$\deg(6) = 3$.  

To obtain the precise value of the finitistic dimension, we further observe: 
The arrow $\beta_7$ ends in the vertex $e_8$ with maximal finite length
$c(e_8) = 7$ of departing paths, and hence
$\pdim (\la e_7 / \la \beta_7) = 1 + 1$-$\deg(7) = 4$.  Consequently,
$\Findim \la = \findim \la = 4$.  \qed
\enddefinition

\head 3.  Generic behavior of the homological dimensions \endhead

Recall that, for any finite dimensional algebra $\Delta$ and $d \in \NN$, the
following affine variety $\modladelta$ parametrizes the $d$-dimensional
$\Delta$-modules:  Let $a_1, \dots, a_r$ be a set of algebra generators for
$\Delta$ over $K$.  For instance, if $\Delta$ is a path algebra modulo
relations, then the primitive idempotents (alias vertices of the quiver),
together with the (residue classes in $\Delta$ of the) arrows constitute
such a set of generators.  For $d \in \NN$, 
$$\modladelta = \{(x_i) \in \prod_{1 \le i \le r} \End_K(K^d)
\mid \text{\ the\ } x_i \text{\ satisfy all relations satisfied by the
\ } a_i\}.$$
  As is well-known, the isomorphism classes of $d$-dimensional (left)
$\la$-modules are in one-to-one correspondence with the orbits of
$\modlad$ under the $\GL_d$-conjugation action.  Indeed, the orbits coincide
with the fibres of the map from $\modladelta$ to the set of isomorphism
classes of $d$-dimensinal left $\Delta$-modules, which maps a point $x$ to
the class of $K^d$, endowed with the $\Delta$-multiplication
$a_i v = x_i(v)$.  If $\C$ is a subvariety of
$\modladelta$, we refer to the modules represented by the points in $\C$ as
{\it the modules in $\C$\/}.

It is, moreover, a standard fact that the homological dimensions of the
$d$-dimensional modules, such as $\pdim$, are generically constant on any
irreducible component of $\modladelta$ (for a proof, see \cite{\CBS, Lemma
4.3} or \cite{\Sma, Theorem 5.3}, where the result is attributed to
Bongartz).  In fact, it is known that, given any irreducible subvariety
$\C$ of $\modladelta$, there exists a dense open subset $U \subseteq \C$ such
that the function $\pdim$ is constant on $U$.   Moreover, this {\it generic
projective dimension on $\C$\/} is the minimum of the projective dimensions
attained on the modules in $\C$.  In most interesting cases, the projective
dimension fails to be constant on all of $\C$, however. (Think, e.g., of the
path algebra $\Delta$ of the quiver $1 \rightarrow 2$, and let
$\C$ be the irreducible component of $\mod_2(\Delta)$, whose points
correspond to the modules with composition factors $S_1, S_2$; here the
generic projective dimension is $0$, while $\pdim(S_1 \oplus S_2) = 1$.) 
This raises the question of how the following generic variant of the
finitistic dimension relates to the classical little finitistic dimension of
$\Delta$.

\definition{Definition 10} The {\it generic left finitistic dimension\/} of a
finite dimensional algebra $\Delta$ is the supremum of the {\rm finite}
numbers gen-$\pdim (\C)$, where $\C$ traces the irreducible components of
the varieties $\modladelta$; here gen-$\pdim (\C)$ is the generic value of
the function $\pdim$, restricted to the modules in $\C$.
\enddefinition 

Clearly, the (left) generic finitistic dimension of an algebra
$\Delta$ is always bounded above by $\findim \Delta$.  When are the two
dimensions equal?  Given an irreducible component $\C \subseteq
\modladelta$, what is the spectrum of values attained by the projective
dimension on $\C$?   

The completeness with which these questions can be answered in the case of a
truncated path algebra $\la$ came as a surprise to us.  The resulting
picture underscores the pivotal role played by tree modules and supplements
the fact that, in the truncated scenario, the irreducible components are
fairly well understood.  They are in one-to-one correspondence with certain
sequences of semisimple modules, as follows:  

Recall that, given a finitely generated left $\la$-module $M$, its {\it
radical layering\/} is $\SS(M) = (J^l M/ J^{l+1} M)_{0 \le l \le L}$.  We
will identify isomorphic semisimple modules so that the radical layerings of
isomorphic $\la$-modules become identical.  That the $K$-dimension of $M$ be
$d$, evidently translates into the equality $\sum_{0 \le l \le L} \dim_K J^l
M/ J^{l+1} M = d$.  For each sequence $\SS = (\SS_0, \dots, \SS_L)$ of
semisimple modules $\SS_l$ with total dimension $d$, let
$\laySS$ be the subset of $\modlad$ consisting of those points which
correspond to the modules with radical layering $\SS$.  Then the locally
closed subvariety $\laySS$ of $\modlad$ is irreducible by
\cite{\BHT, Theorem 5.3}, whence so is its closure in $\modlad$.  The
maximal candidates among the closures $\overline{\laySS}$, where $\SS$
traces the sequences $\SS$ of total dimension $d$, are therefore the
irreducible components of
$\modlad$; indeed, there are only finitely many such sequences.  It is,
moreover, easy to recognize whether a given sequence
$\SS$ of semisimple modules as above arises as the radical layering of a
$\la$-module, that is, whether $\laySS \ne \varnothing$ (see \cite{\BHT}). 
Namely, suppose that $\SS_l = \bigoplus_{0 \le l \le L} S_i^{(i,l)}$ and let
$P$ be the projective cover of $\SS_0$.  Then
$\laySS \ne \varnothing$ if and only if there exists a set $\S$ of paths in
$P$, which is closed under initial subpaths, such that $\S$ is {\it
compatible with $\SS$\/} in the following sense:  For each $i \in \{1,
\dots, n\}$ and each $l \in \{0, 1, \dots, L\}$, the set $\S$ contains
precisely $s(i,l)$ paths of length $l$ which end in the vertex
$e_i$.  Observe that, whenever $M$ is a module with radical layering $\SS(M)
= \SS$, any skeleton of $M$ is compatible with $\SS$.  Consequently, the
requirement that $\laySS \ne \varnothing$ implies that the $\l$-th layer
$\SS_l$ of $\SS$ be a direct summand of the $l$-th layer $J^l P/ J^{l+1} P$
in the radical layering of $P$.

\proclaim{Theorem 11} Let $\SS = (\SS_0, \SS_1, \dots, \SS_L)$ be a sequence
of semisimple $\la$-modules such that $\laySS \ne \varnothing$, and let $P$ a
projective cover of $\SS_0$.  Moreover, suppose 
$$J^l P / J^{l+1}P = \biggl (\bigoplus_{1 \le i \le n} S_i^{s(i,l)}
\biggr) \oplus \biggl( \bigoplus_{1 \le i \le n} S_i^{r(i,l)} \biggr)$$  for
suitable nonnegative integers
$r(i,l)$; here $s(i,l)$ is the multiplicity of $S_i$ in $\SS_l$ as above.  
\smallskip  

\noindent {\rm (1)}  The projective dimension of a module
$M$ depends only on its radical layering $\SS(M)$. In other words, the
projective dimension is constant on each of the varieties $\laySS$.    This
constant value, denoted $\pdim \SS$, is the generic projective dimension of
the irreducible subvariety $\overline{\laySS}$ of
$\modlad$.
\smallskip  

\noindent {\rm (2)} If $\pdim \SS > 0$, then 
$$\pdim \SS = 1 \, +\, \sup\{l\text{-}\deg(c(e_i)) \mid i \le n,\, l \le L
\text{\ with \ } r(i,l) \ne 0\}.$$  {\rm (}We adopt the standard convention
``$1 + \infty =
\infty$".{\rm )} In particular, $\pdim \SS$ is finite if and only if
$r(i,l) = 0$ for all cyclebound vertices $e_i$, that is, if and only if every
simple module of infinite projective dimension has the same composition
multiplicity in  
$P$ as in $\bigoplus_{0 \le l \le L} \SS_l$.   
\smallskip

\noindent {\rm (3)}  The generic finitistic dimension of $\la$ coincides
with $\, \findim \la$.  It is the projective dimension of a tree module
$\T$  --  of dimension $d$ say  --   whose orbit closure is an irreducible
component of $\modlad$. 
\endproclaim

Computing $\pdim \SS$ in concrete examples amounts to performing at most $n$
counts:  Indeed, if $r(i,l) \ne 0$ for some $l$, then $l$-$\deg(c(e_i)) \le
l_i$-$\deg(c(e_i))$, where $l_i$ is maximal with $r(i, l_i) \ne 0$.  Observe
moreover that the event $\pdim \SS = 0$ is readily recognized: It occurs if
and only if $\SS = \SS(P)$; in this case,
$\laySS$ consists of the $\GL_d$-orbit of $P$ only.   

We smooth the road towards a proof of Theorem 11 with two preliminary
observations.     

\proclaim{Observation 12}  Given any finitely generated $\la$-module with
skeleton $\S$, there exists a direct sum of tree modules with the same
skeleton.

In particular, the syzygy of any finitely generated $\la$-module is
isomorphic to the syzygy of a direct sum of tree modules, and all projective
dimensions in  $\{0, 1, \dots, \findim \la\}$  are attained on tree modules.
\endproclaim 
   
\demo{Proof}  Let $M$ be any finitely generated left $\la$-module, $P =
\bigoplus_{1 \le r \le t} \la z_r$ a projective cover of $M$ with $z_r = e(r)
\in
\{e_1, \dots, e_n\}$, and $\sigma \subseteq P$ a skeleton of $M$.  For fixed
$r \le t$, let $\sigma^{(r)}$ be the subset of $\sigma$ consisting of all
paths in $\S$ of the form $pz_r$.  Then $\T^{(r)} := \la z_r /(\sum_{q\
\sigma^{(r)}\text{-critical}} \la q)$ is a tree module whose branches are
precisely the maximal paths in $\sigma^{(r)}$ relative to the ``initial
subpath order".  Hence $\bigoplus_{1 \le r \le t} \T^{(r)}$ is a direct sum
of tree modules, again having skeleton $\sigma$.  Since, by Theorem 1, any
skeleton of a module determines its syzygy up to isomorphism, the remaining
claims follow. \qed
\enddemo      

The next observation singles out  candidates for the tree module postulated
in Theorem 11(3).   Let $\epsilon$ be the sum of all non-cyclebound
primitive idempotents in the full set $e_1, \dots, e_n$.  (In Example 3, we
have $\epsilon = e_8 + \dots + e_{15}$.)  Clearly, the left ideal $\la
\epsilon
\subseteq \la$ of finite projective dimension equals $\epsilon \la
\epsilon$.  In particular, given any left
$\la$-module $M$, the subspace
$\epsilon M$ is a sub{\it module\/} of $M$.

\proclaim{Observation 13} Let $e_i$ be any vertex of $Q$. Then $\pdim
\epsilon J e_i < \infty$, and
$$\pdim \epsilon J e_i \ge \pdim \la q,$$ for every nonzero path $q$ of
positive length in $\la$ which starts in $e_i$ and satisfies $\pdim \la q <
\infty$.

Moreover: The factor module $\T_i = \la e_i / \epsilon J e_i$ is a tree
module.  If $\dim_K \T_i = d_i$, and $\SS(\T_i) =
\SS^{(i)}$ is the radical layering of $\T_i$, then the subvariety $\mod
\bigl(\SS^{(i)}\bigr)$ of
$\mod_{d_i}(\la)$ coincides with the $\GL_{d_i}$-orbit of 
$\T_i$ and is open in $\mod_{d_i}(\la)$.
\endproclaim

\demo{Proof}  We first address the second set of claims.  Let $p_{ij}$, $1
\le j \le t_i$, be the different paths of positive length in
$\la$ which start in
$e_i$, end in a non-cyclebound vertex, and are minimal with these properties
in the ``initial subpath order"; that is, every proper intial subpath of
positive length of one of the $p_{ij}$ ends in a cyclebound vertex.  Clearly,
$\epsilon J e_i = \bigoplus_{1 \le j \le t_i} \la p_{ij}$, which shows in
particular that $\T_i$ is a tree module.  Moreover, any module $M$ sharing
the radical layering of $\T_i$ also has projective cover $\la e_i$, and a
comparison of composition factors\ shows that every epimorphism $\la e_i
\rightarrow M$ has kernel $\epsilon J e_i$.  Thus $M \cong T_i$, which shows
$\mod \bigl(\SS^{(i)}\bigr)$ to equal the
$\GL_{d_i}$-orbit of $\T_i$.  Moreover, it is readily checked that
$\Ext^1_\la(\T_i, \T_i) = 0$, whence the orbit 
$\mod \bigl(\SS^{(i)}\bigr)$ of $\T_i$ is open in $\mod_{d_i}$ (see
\cite{\Del, Corollary 3}), and the proof of the final assertions is
complete.          

For the first claim, let $q = q e_i$ be a nonzero path of positive length in
$\la$ with $\pdim \la q < \infty$.  Then $q$ ends in a non-cyclebound vertex
by Lemma 8  --  call it $e$  --  and hence $q$ has an initial subpath $q'$
among the paths $p_{ij}$; let $e'$ be the (non-cyclebound) terminal vertex
of $q'$.  If $l$ and $l'$ are the lengths of $q$ and $q'$, respectively, 
$c(e') - c(e) \ge l - l' \ge 0$, and hence $c(e') + l' \ge c(e) + l$.  This
shows
$$\pdim \la q' = l'\text{-}\deg(c(e')) \ge l\text{-}\deg(c(e))  = \pdim \la
q,$$  which yields the desired inequality.  \qed \enddemo

\demo{Proof of Theorem 11}  (1)  Let $M$ be a module with radical layering
$\SS$ and $\sigma$ any skeleton of $M$.  By \cite{\BHT, Theorem 5.3}, the
points in $\modlad$ parametrizing the modules that share this skeleton
constitute a dense open subset of $\laySS$.  All modules represented by this
open subvariety have the same projective dimension as
$M$, because any skeleton of a module pins down its syzygy up to
isomorphism.  Therefore,
$\pdim M$ is the generic value of the function $\pdim$ on the irreducible
subvariety $\overline{\laySS}$ of $\modlad$. 

(2) Suppose that $\pdim \SS > 0$, which means $r(i,l) > 0$ for some pair
$(i,l)$.  Let $M$ be any module with $\SS(M) = \SS$.  By part (1), $\pdim
\SS = \pdim M$.  To scrutinize the projective dimension of $M$, let 
$\widehat{\S}$ be a skeleton of $P$ and $\S \subset \widehat{\S}$ a skeleton
of $M$.  We have
$\Omega^1(M) \cong \bigoplus_{q\ \S\text{-}\text{critical}} \la q$ by
Theorem 2.  Since $r(i,l) > 0$ whenever $q$ is a $\S$-critical path of length
$l$ ending in $e_i$, we glean that $\pdim M$ is bounded above by the
supremum displayed in part (2) of Theorem 11.  For the reverse inequality,
choose any pair $(i,l)$ with $r(i,l) > 0$.  This inequality amounts to the
existence of a path $p z_r$ of length $l$ in $\widehat{\S} \setminus
\S$ which ends in $e_i$.  Denote by
$p' z_r$ the maximal initial subpath of $p z_r$ which belongs to $\S$. Since
$p z_r \notin \S$, there is a unique arrow $\alpha$ such that
$\alpha p' z_r$ is in turn an initial subpath of $p z_r$.  In particular, if
$q = \alpha p'$, then $q z_r$ is a $\S$-critical path ending in some vertex
$e_j$.  Invoking once again the above decomposition of
$\Omega^1(M)$, we deduce that the cyclic left ideal $\la q$ is isomorphic to
a direct summand of $\Omega^1(M)$.  By Theorem 5, it therefore suffices to
show that  the $\len(q)$-degree of $c(e_j)$ is larger than or equal to
$l$-$\deg(c(e_i))$.   For that purpose, we write $p z_r = q' q z_r$ for a
suitable path $q'$ in $\la$.  Since $c(e_j) \ge c(e_i) + \len(q')$, we obtain
$c(e_j) \ge c(e_i)$, and consequently $c(e_j) + \len(q) \ge c(e_i) +  l$.  We
conclude 
$$\left[\frac{c(e_j)}{L+1}
\right] + \left[\frac{c(e_j) + \len(q)}{L+1} \right]  \ge
\left[\frac{c(e_i )}{L+1} \right] + \left[\frac{c(e_i) + l}{L+1} \right] =
l\text{-}\deg(c(e_i)).$$  Thus $\pdim M - 1 \ge l$-$\deg(c(e_i))$ as
required.  The final equivalence under (2) is an immediate consequence.

(3)  By construction, the tree modules $\T_i$ of Observation 13 all have
finite projective dimension.  Combining the first part of this observation
with the final statement of Theorem 2, we moreover see that $\findim \la$
equals the maximum of these dimensions.  The final statement of Observation
13 now completes the proof of (3).  \qed         
\enddemo

Let $\SS = (\SS_0, \dots, \SS_l)$ again be a sequence of semisimple modules
of total dimension $d$ such that $\laySS \ne \varnothing$.  As we saw, the
projective dimension $\pdim \SS$ holds some  information about path lengths
in $KQ$; namely on the lengths of paths starting in vertices that belong to
the support of
$\Omega^1(M)$, where $M$ is any module in $\lay\SS$.  To obtain a tighter
correlation between $Q$ and the homology of $\la$, we will next explore the
full spectrum of values of the function
$\pdim$ attained on the closure $\overline{\laySS}$.  While those ranges of
values are better gauges of how the vertices corresponding to the simples in
the various layers $\SS_l$ of $\SS$ are placed in the quiver $Q$, the
refined homological data still do not account for the intricacy of the
embedding of $\overline{\laySS}$ into $\modlad$ in general.  (See the
comments following the next theorem.)  On the other hand, for $\pdim \SS <
\infty$ and small
$L$, far more of this picture is preserved in the homology than in the
hereditary case.      

We first recall from \cite{\BHT, Section 2.B} that, for any $M$ in
$\overline{\laySS}$, the sequence $\SS(M)$ is larger than or equal to $\SS$
in the following partial order:  Suppose that $\SS$ and
$\SS'$ are semisimple modules with $\bigoplus_{0 \le l \le L} \SS_l = 
\bigoplus_{0 \le l \le L} \SS'_l$.  Then ``$\SS' \ge \SS$" means that
$\bigoplus_{l \le r} \SS_l$ is a direct summand of $\bigoplus_{l \le r}
\SS'_l$, for all $r \ge 0$.  In intuitive terms this says that, in the
passage from
$\SS$ to $\SS'$, the simple summands of the $\SS_l$ are only upwardly mobile
relative to the layering. 

\proclaim{Lemma 14} If $\SS' \ge \SS$ and $\mod(\SS') \ne \varnothing$, then
$\pdim \SS' \ge \pdim \SS$. \endproclaim 

\demo{Proof} Let $P$ be a projective cover of $\SS_0$ as before and
$P'$ a projective cover of $\SS'_0$.  Decompose the radical layers of $P'$ in
analogy with the decomposition given for $P$ above:  
$$J^lP' / J^{l+1} P' = \bigoplus_{1 \le i \le n} S_i^{s'(i,l)} \oplus 
\bigoplus_{1 \le i \le n} S_i^{r'(i,l)},$$ where $\SS'_l = \bigoplus_{1 \le
i \le n} S_i^{s'(i,l)}$.  It follows immediately from the definition of the
partial order of sequences of semisimples that, whenever $r(i,l) > 0$, there
exists $l' \ge l$ with
$r'(i,l') > 0$.  In light of Theorem 11, this proves the lemma. \qed
\enddemo   

We will give two descriptions of the range of values of $\pdim$ on the
closure $\overline{\laySS}$.  For a combinatorial version, we keep the
notation of Theorem 11 and the proof of Lemma 14:  Namely, $\SS_l =
\bigoplus_{1 \le i \le n} S_i^{s(i,l)}$, and
$P$ is a projective cover of $\SS_0$.  From $\laySS
\ne \varnothing$, one then obtains $J^l P/ J^{l+1} P = \SS_l \oplus
\bigoplus_{1 \le i \le n} S_i^{r(i,l)}$.   In our graph-based description of
the values $\pdim M > \pdim \SS$, where $M$ traces
$\overline{\laySS}$, the exponents
$s(i,l)$ take over the role played by the $r(i,l)$ relative to the generic
projective dimension, $\pdim \SS$:  Recall from Theorem 11 that, whenever
$\pdim \SS$ is nonzero, it is the maximum of the values $1 +
l$-$\deg(c(e_i)) \in \NN \cup \{0, \infty\}$ which accompany the pairs
$(i,l)$ with $r(i,l) > 0$.  (Note: In view of $\SS_0 = P/JP$, the inequality
$r(i,l) > 0$ entails $l\ge 1$.) 

Now, we consider the different candidates $n_1,\dots, n_v$ among those
elements in $\NN \cup \{0, \infty
\}$, which have the form 
$$1 + l\text{-}\deg(c(e_j)), \  l \ge 1,  \ S_j
\subseteq \SS_l$$ and are {\it strictly larger\/} than $\pdim \SS$. In other
words, 
$$\{n_1, \dots, n_v\} = [\pdim \SS + 1\,,\ \infty] \cap \{1 + 
l\text{-}\deg(c(e_j)) \mid l \ge 1, s(j,l) > 0 \}.$$

\proclaim{Theorem 15}  Let $\SS$ be a semisimple sequence of total dimension
$d$ with $\laySS \ne \varnothing$.  
\smallskip

\noindent The range of values,
$$\{ \pdim M \mid M \text{\ in \ } \overline{\laySS} \},$$ of the function
$\pdim$ on the closure of $\laySS$ in $\mod_d$ is equal to the following
coinciding sets:

$$\{\pdim \SS' \mid \SS' \ge \SS,\, \mod(\SS') \ne
\varnothing \} = \{\pdim \SS\} \cup \{n_1, \dots, n_v\}.$$ 
\smallskip      
\endproclaim

In general, describing the closure of
$\laySS$ in $\modlad$ is an intricate representation-theoretic task, a fact
not reflected by the homology.  For instance: 
$\bullet$  When $\SS'$ is a sequence of semisimple modules such that $\SS'
\ge \SS$ and $\mod(\SS') \ne \varnothing$, the intersection
$\overline{\laySS} \cap
\mod(\SS')$ may still be empty.  $\bullet$ The condition $\overline{\laySS}
\cap \mod(\SS') \ne \varnothing$ does not imply $\mod(\SS')
\subseteq \overline{\laySS}$.  See the final discussion of our example for
illustration.   

\demo{Proof}  Set $\P = \{ \pdim M \mid M \text{\ in \ } \overline{\laySS}
\}$.  We already know that $\P \subseteq \{ \pdim \SS'
\mid \SS' \ge \SS\}$; indeed, this is immediate from Lemma 14 and the
remarks preceding it.  

Suppose that $\SS'$ is a sequence of semisimple modules with $\SS' \ge \SS$
and $\mod(\SS') \ne \varnothing$.  Assume $\pdim \SS' > \pdim \SS$, which, in
particular, implies $\pdim \SS' > 0$. To show that $\pdim \SS'$ equals one
of the $n_k$, we adopt the notation used in the proof of Lemma 14.  By
Theorem 11,
$\pdim \SS' = 1 + a$-$\deg(c(e_i))$ for some pair $(i,a)$ with $r'(i,a) >
0$.  Again invoking Theorem 11, we moreover infer that $r(i,a) = 0$ from
$\pdim \SS < \pdim \SS'$.  If $s(i,a) > 0$, we are done, since necessarily
$a \ge 1$.  So let us suppose that also $s(i,a) = 0$, meaning that $S_i$
fails to be a summand of the $a$-th layer $J^a P / J^{a+1} P$ of
$P$.  In light of $S_i \subseteq J^a P'/ J^{a+1} P'$, this entails the
existence of a simple $S_j \subseteq \SS'_0 / \SS_0$ with the property that
$S_i \subseteq J^a e_j / J^{a+1} e_j$.  Consequently,
$c(e_j) \ge c(e_i) + a$.  On the other hand, $S_j \subseteq \bigoplus_{l \ge
1} \SS_l$, because the total multiplicities of the simple summands of
$\SS$ and $\SS'$ coincide.  This means $s(j,k) > 0$ for some $k \ge 1$.  In
light of $\SS_0 \oplus S_j \subseteq \SS'_0$ and $\bigoplus_{0 \le l \le L}
\SS_l = \bigoplus_{0 \le l \le L} \SS'_l$, we deduce that $r'(j,b) > 0$ for
some pair $(j,b)$ with $b \ge 1$ and $s(j,b) > 0$.  Another application of
Theorem 11 thus yields $\pdim \SS' - 1 \ge b\text{-}\deg(c(e_j))$ $=$
$\left[\frac{c(e_j)}{L+1} \right] +
\left[\frac{c(e_j) + b}{L+1} \right]$ $\ge$ $\left[\frac{c(e_i) + a}{L+1}
\right] + \left[\frac{c(e_i) + a + b}{L+1} \right] \ge a\text{-}\deg(c(e_i))$
$=$ $\pdim \SS' - 1$.  We conclude that all inequalities along this string
are actually equalities, that is, $b\text{-}\deg(c(e_j)) =
a\text{-}\deg(c(e_i))$.  This shows that $\pdim \SS' = 1 +
b\text{-}\deg(c(e_j))$ for a pair $(j,b)$ with $s(j,b) > 0$ as required.

Finally, we verify that each of the numbers $n_k$ belongs to $\P$. By
definition, $n_k$ is of the form $1 + l$-$\deg(c(e_i))$ for some pair
$(i,l)$ with $l \ge 1$ and  $s(i,l) > 0$.  Let
$D$ be any direct sum of tree modules with $\SS(D) = \SS$; in light of
$\laySS \ne \varnothing$, such a module $D$ exists by Observation 12.  Then
there is a tree direct summand $\T$ of $D$ with a branch that contains an
initial subpath $q$ of length $l$ ending in the vertex $e_i$.  The direct
sum of tree modules $D' = (\T/ \la q) \oplus \la q \oplus D/\T$ belongs to
$\overline{\laySS}$.  In fact, $D'$ is well known to belong to the closure
of the orbit of $D$ in $\modlad$; see, e.g., \cite{\CB, Section 3, Lemma 2}.
Therefore $\pdim D' \in \P$.  As for the value  of this projective
dimension: Up to isomorphism,
$\la q$ is a direct summand of the syzygy of the tree module $\T/ \la q$: 
indeed, $q$ is $\S$-critical relative to the obvious skeleton $\S$ of $\T/
\la q$ consisting of all initial subpaths of the branches (see Theorem 2 and
the comments accompanying Definition 6).  Theorems 5 and 11 moreover yield
$$ \pdim D' - 1 \ge \pdim \la q = l\text{-}\deg (c(e_i)) = n_k - 1 > \pdim
\SS - 1 =\pdim D - 1,$$  whence $\pdim(\T/ \la q) = n_k = \pdim D'$.  This
shows $n_k$ to belong to $\P$ and completes the argument.  \qed             
\enddemo
 
\definition{A final visit to Example 3}  (a)  First, let $\SS = \SS(\la
e_1)$ be the radical layering of the projective tree module $\T = \la e_1$. 
By Theorems 11 and 15, the values of $\pdim$ on
$\overline{\laySS}$ are $\pdim \SS = 0, 2, 3, 4, \infty$.  For instance,
$\pdim \bigl( (\T/ \la \beta_3 \alpha_2 \beta_1) \oplus (\la \beta_3 \alpha_2
\beta_1) \bigr)$ equals $2$.  Note that the value $1$, on the other hand, is
not attained. 
\smallskip

\noindent (b)  Next we justify the comments following Theorem 14.  Let $\SS$
and $\SS'$ be the radical layerings of the modules $M$ and
$M'$ with the following graphs, respectively:  

$$\xymatrixrowsep{0.6pc} \xymatrixcolsep{0.5pc} 
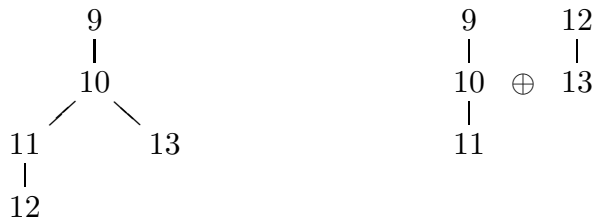
\xymatrix{
& 9 \ar@{-}[d] \\
& 10 \ar@{-}[dl] \ar@{-}[dr] \\ 
11 \ar@{-}[d] & & 13 \\ 12
} \qquad\qquad \qquad\qquad
\xymatrix{
9 \ar@{-}[d] & & 12 \ar@{-}[d] \\ 
10 \ar@{-}[d] &\save+<0ex,0ex> \drop{\oplus} \restore \restore & 13 \\ 
11
}$$

\noindent Then $\SS' \ge \SS$, while
$\overline{\laySS} \cap \mod(\SS') = \varnothing$.  On the other hand, if
$M$, $M'$, and $M''$ are given by the graphs

$$\xymatrixrowsep{0.6pc} \xymatrixcolsep{0.5pc} 
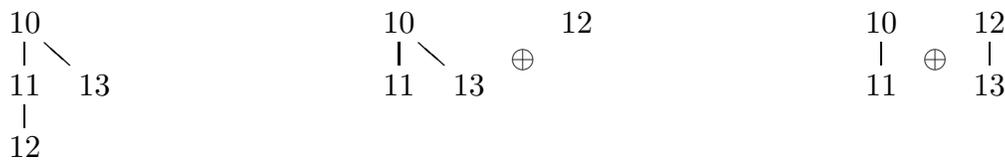
\xymatrix{ 
10 \ar@{-}[d]
\ar@{-}[dr] \\ 
11 \ar@{-}[d]  & 13 \\ 12
} \qquad\qquad \qquad\qquad 
 \xymatrix{
10 \ar@{-}[d] \ar@{-}[dr] && \save+<0ex,-2.5ex> \drop{\oplus} \restore & 12 \\ 
11 & 13
} \qquad\qquad \qquad\qquad 
 \xymatrix{
10 \ar@{-}[d] & \save+<0ex,-2.5ex> \drop{\oplus} \restore & 12 \ar@{-}[d] \\ 
11 & & 13
} $$
\noindent then $\SS' := \SS(M')$ equals $\SS(M'')$, and the intersection
$\overline{\laySS}
\cap \mod(\SS')$ contains $M'$, but not $M''$.
\enddefinition


\Refs
\widestnumber\key{{\bf 99}}

\ref\no \BHT \by E. Babson, B. Huisgen-Zimmermann, and R. Thomas \paper
Generic representation theory of quivers with relations
\paperinfo manuscript 2007. 
\endref

\ref\no \But \by M. Butler \paper The syzygy theorem for monomial algebras
\inbook Trends in the Representation Theory of Finite Dimensional Algebras
(Seattle, Washington 1997) \eds E.L. Green and B. Huisgen-Zimmermann
\bookinfo Contemp. Math. \vol 229   \publ Amer. Math. Soc.  \publaddr
Providence \yr 1998  \pages 111-116
\endref

\ref\no \CB  \by W. Crawley-Boevey \paper Lectures on representations of
quivers
\inbook Course Notes 1992  \paperinfo
http://www.maths.leeds.ac.uk/$\sim$pmtwc/  
\endref

\ref\no \CBS  \by W. Crawley-Boevey and J. Schr\"oer \paper Irreducible
components of varieties of modules \jour J. reine angew. Math. \vol 553 \yr
2002 \pages 201-220  \endref

\ref\no \Del \by Jose Antonio de la Pe\~na \paper  Tame algebras: some
fundamental notions \paperinfo Erg\"anzungsreihe 95-010
Sonderforschungsbereich 343, Universit\"at Bielefeld \yr 1995 \endref 

\ref\no \GKK \by E.L. Green, E.E. Kirkman, and J.J. Kuzmanovich \paper
Finitistic dimensions of finite dimensional monomial algebras \jour J.
Algebra \vol 136 \yr 1991 \pages 37-51 \endref  

\ref\no \pre \by B. Huisgen-Zimmermann \paper Predicting syzygies over
monomial relation algebras \jour Manuscripta Math. \vol 70 \yr 1991 \pages
157-182 \endref 

\ref\no\dom \bysame \paper Homological domino effects and the first
finitistic dimension conjecture \jour Invent. Math. \vol 108 \yr 1992 \pages
369-383 
\endref

\ref\no \Sma \by S. O. Smal\/o \paper Lectures on Algebras, Mar del Plata,
Argentina, March 2006  \paperinfo Revista de la Uni\'on Matem\'atica
Argentina, to appear \endref


\endRefs
\enddocument